\documentclass{article}

\usepackage{amsmath,amssymb,amsfonts,amsthm}
\usepackage[all]{xy}

\newtheorem{defn}{Definition}[section]

\newtheorem{prop}[defn]{Proposition}
\newtheorem{cor}[defn]{Corollary}
\newtheorem{rem}[defn]{Remark}
\newtheorem{thm}[defn]{Theorem}
\newtheorem{lemma}[defn]{Lemma}

\newcommand{\longto}[1][]{\stackrel{#1}{\longrightarrow}}
\newcommand{\otgnol}[1][]{\stackrel{#1}{\longleftarrow}}
\newcommand{\integers}{\mathbb{Z}}
\newcommand{\calA}{\mathcal{A}}
\newcommand{\calE}{\mathcal{E}}
\newcommand{\calO}{\mathcal{O}}
\newcommand{\calM}{\mathcal{M}}

\newcommand{\coarseM}{\mathrm{M}}
\newcommand{\rk}{\mathrm{rk}}
\newcommand{\CH}{\mathrm{CH}}
\newcommand{\NS}{\mathrm{NS}}
\newcommand{\SL}{\mathrm{SL}}
\newcommand{\GL}{\mathrm{GL}}
\newcommand{\PSL}{\mathrm{PSL}}
\newcommand{\PGL}{\mathrm{PGL}}
\newcommand{\Gm}{\mathbb{G}_m}
\newcommand{\Quot}{\mathrm{Quot}}
\newcommand{\Drap}{\mathrm{Drap}}
\newcommand{\characteristic}{\mathrm{char}}
\newcommand{\Schemes}{\underline{Sch}}
\newcommand{\Sets}{\underline{Sets}}
\newcommand{\Spec}{\mathrm{Spec}\:}
\renewcommand{\H}{\mathrm{H}}
\newcommand{\Endsheaf}{\mathit{End}}
\newcommand{\Homsheaf}{\mathit{Hom}}
\newcommand{\End}{\mathrm{End}}
\newcommand{\Hom}{\mathrm{Hom}}
\newcommand{\Ext}{\mathrm{Ext}}
\newcommand{\Pic}{\mathrm{Pic}}
\newcommand{\m}{\mathfrak{m}}
\newcommand{\obstr}{\mathrm{ob}}
\newcommand{\supp}{\mathrm{supp}}

\newcommand{\ks}{\mathrm{ks}}
\newcommand{\op}{\mathrm{op}}
\newcommand{\dual}{\mathrm{dual}}
\newcommand{\tr}{\mathrm{tr}}
\newcommand{\lp}{\mathrm{lp}}

\newcommand{\Mat}{\mathrm{Mat}}

\begin{document}  \bibliographystyle{plain}

\title{Moduli schemes of generically simple\\Azumaya modules}
\author{Norbert Hoffmann\thanks{Mathematisches Institut der Universit\"at G\"ottingen, Bunsenstr. 3--5, D-37073 G\"ottingen, Germany. email: hoffmann@uni-math.gwdg.de}
     \and Ulrich Stuhler\thanks{Mathematisches Institut der Universit\"at G\"ottingen, Bunsenstr. 3--5, D-37073 G\"ottingen, Germany. email:  stuhler@uni-math.gwdg.de}
}
\date{}
\maketitle

\section*{Introduction}

Let $X$ be a smooth projective variety, e.\,g.\ a surface, over an algebraically closed field $k$. Let $\calA$ be a sheaf of Azumaya algebras over $X$ or more generally, a torsion free
coherent sheaf of algebras over $X$ whose generic fiber $\calA_{\eta}$ is a central simple algebra over the function field of $X$. This paper is about moduli schemes of generically simple,
locally projective $\calA$-module sheaves $E$.

These moduli schemes are in close analogy to the Picard variety of $X$. In fact, our main result says that we do not need any stability condition for our sheaves $E$ to construct coarse moduli
schemes parameterizing them, say with fixed Hilbert polynomial or Chern classes. We find that these schemes are in general not proper over $k$, but they have natural compactifications:
Working with torsion free sheaves $E$ instead of only locally projective ones, we obtain projective moduli schemes.

This gives lots of interesting moduli spaces, which certainly deserve further study. For example, we show in section \ref{deformations} that they are smooth projective and even symplectic
if $X$ is an abelian or K3 surface and $\calA$ is an Azumaya algebra. They are also related to classifying isomorphism types of Azumaya algebras $\calA$ in a given central division algebra
$\calA_{\eta} = D$, a topic already present in the classical literature on algebras; this relation is explained in section \ref{orders}.

We construct these moduli schemes in section \ref{construction}. We use standard techniques from geometric invariant theory (GIT) and a boundedness result, which has been known for some time
in the case of characteristic $\characteristic( k) = 0$, but is one of the deep results of A. Langer in \cite{langer} for $\characteristic( k) = p > 0$. Our construction works for any integral
projective scheme $X$ over $k$; the precise setup is formulated in section \ref{modules}. As M. Lieblich has pointed out to us, this construction can also be seen as a special case of
Simpson's general result \cite[Theorem 4.7]{simpson}, at least if $\characteristic( k) = 0$.

In section \ref{deformations}, we apply deformation theory to our $\calA$-modules, mainly in the case where $X$ is a surface and $\calA$ is an Azumaya algebra. Besides the smoothness
mentioned above, we also show here that torsion free sheaves are really necessary to obtain projective moduli schemes, because locally projective sheaves of $\calA$-modules can degenerate
to torsion free ones.

During the final preparations of this paper, we were informed about the MIT-thesis of M. Lieblich \cite{lieblich1, lieblich2}. This thesis works much more systematically and abstractly and
contains several results similar to ours in the language of algebraic stacks. We thank A. J. de Jong and M. Lieblich for informations concerning their work. Similar results have also been
obtained independently by K. Yoshioka \cite{yoshioka}; they have been used by D. Huybrechts and P. Stellari \cite{huybrechts-stellari} to prove a conjecture of Caldararu.

\section{Families of $\calA$-modules} \label{modules}

Let $X$ be an integral projective scheme over the algebraically closed field $k$. Throughout this paper, $\calA$ denotes a sheaf of associative $\calO_X$-algebras satisfying the
following properties:
\begin{itemize}
 \item[1.] As a sheaf of $\calO_X$-modules, $\calA$ is coherent and torsion free.
 \item[2.] The stalk $\calA_{\eta}$ of $\calA$ over the generic point $\eta \in X$ is a central simple algebra over the function field $F = k( X) = \calO_{X, \eta}$.
\end{itemize}
For example, $X$ could be a smooth projective variety over $k$, and $\calA$ could be a sheaf of Azumaya algebras over $X$.
\begin{rem} \upshape
  If $\dim X = 1$, then $\calA_{\eta}$ is a matrix algebra over $k( X)$ by Tsen's theorem. So the first interesting case is $\dim X = 2$.
\end{rem}

Our main objects will be generically simple torsion free $\calA$-modules, i.\,e.\ sheaves $E$ of left $\calA$-modules over $X$ which are torsion free and coherent as $\calO_X$-modules and
whose generic fiber $E_{\eta}$ is a simple module over the central simple algebra $\calA_{ \eta}$. By Wedderburn's structure theorem, we have $A_{\eta} \cong \Mat( n \times n; D)$ for a
division algebra $D$, say of dimension $r^2$ over $k( X)$; that $E_{\eta}$ is simple means that it is Morita equivalent to a one-dimensional vector space over $D$. In particular, $E$ has rank
$r^2 n$ over $\calO_X$.

Note that any such $\calA$-module $E$ has only scalar endomorphisms: Indeed, $\End_{\calA}( E)$ is a finite-dimensional $k$-algebra; it has no zero-divisors because it embeds into the
division algebra $\End_{\calA_{\eta}}( E_{\eta}) \cong D^{\op}$. This implies $\End_{\calA}( E) = k$ because $k$ is algebraically closed.

\begin{lemma} \label{nonsplit}
  Suppose that $k \subseteq K$ is a field extension, and let $K( X_K)$ be the function field of $X_K := X \times_k \Spec K$. If $D$ is a finite-dimensional division algebra over $k( X)$,
  then $D_K := D \otimes_{k( X)} K( X_K)$ is a division algebra, too.
\end{lemma}
\begin{proof}
  Since $k$ is algebraically closed, $k( X) \otimes_k K$ is an integral domain; its quotient field is $K( X_K)$. Suppose that $D_K$ contains zero divisors. Clearing denominators, we can then
  construct zero divisors in $D \otimes_{k( X)} (k( X) \otimes_k K)$, which is clearly isomorphic to $D \otimes_k K$. Consequently, there is a finitely generated $k$-algebra $A \subseteq K$
  such that $D \otimes_k A$ contains zero divisors. These zero divisors are automatically nonzero modulo some maximal ideal $\m \subset A$, so $D \otimes_k A/\m$ also contains zero divisors.
  But $A/\m \cong k$ by Hilbert's Nullstellensatz; hence $D$ contains zero divisors. This contradiction shows that $D_K$ has to be a division algebra if $D$ is.
\end{proof}

\begin{cor} \label{nonsplit2}
  If $E$ is a generically simple torsion free $\calA$-module, then the pullback $E_K$ of $E$ to $X_K$ is a generically simple torsion free module under the pullback $\calA_K$ of $\calA$.
\end{cor}
\begin{proof}
  $E_K$ is clearly torsion free and coherent over $\calO_{X_K}$. Since the generic fiber of $E$ is Morita equivalent to a one-dimensional $D$-vector space, the generic fiber of $E_K$ is Morita
  equivalent to a one-dimensional $D_K$-vector space; hence $E_K$ is generically simple.
\end{proof}

\begin{defn}
  A \emph{family of generically simple torsion free $\calA$-modules} over a $k$-scheme $S$ is a sheaf $\calE$ of left modules under the pullback $\calA_S$ of $\calA$ to $X \times_k S$
  with the following properties:
  \begin{itemize}
   \item[1.] $\calE$ is coherent over $\calO_{X \times_k S}$ and flat over $S$.
   \item[2.] For every point $s \in S$, the fiber $\calE_s$ is a generically simple torsion free $\calA_{k (s)}$-module.
  \end{itemize}
  Here $k(s)$ is the residue field of $S$ at $s$, and the fiber $\calE_s$ is by definition the pullback of $\calE$ to $X \times_k \Spec k(s)$.
\end{defn}

We denote the corresponding moduli functor by
\begin{displaymath}
  \calM = \calM_{\calA/X}: \Schemes_k \longto \Sets;
\end{displaymath}
it sends a $k$-scheme $S$ to the set of isomorphism classes of families $\calE$ of generically simple torsion free $\calA$-modules over $S$. Our main goal is to construct and study coarse
moduli schemes for this functor.

If $\calE$ is a family of generically simple torsion free $\calA$-modules over $S$, then there is an open subset of $S$ above which these $\calA$-modules are locally projective. However, we
work with all torsion free $\calA$-modules because they satisfy the following valuative criterion for properness:
\begin{prop} \label{proper}
  Let $V$ be a discrete valuation ring over $k$ with quotient field $K$. Then the restriction map
  \begin{displaymath}
    \calM( \Spec V) \longto \calM( \Spec K)
  \end{displaymath}
  is bijective.
\end{prop}
\begin{proof}
  Let $\pi \in V$ be a uniformising element, and let $l = V/(\pi)$ be the residue field of $V$. We denote by
  \begin{displaymath}
    X_K \longto[j] X_V \otgnol[i] X_l
  \end{displaymath}
  the open embedding of the generic fiber and the closed embedding of the special fiber; here $X_A := X \times_k \Spec A$ for any $k$-algebra $A$. Let $\eta$ (resp. $\xi$) be the generic
  point of $X_K$ (resp. of $X_l$), and let
  \begin{displaymath}
    j_{\eta}: \Spec \calO_{X_V, \eta} \longto X_V \quad (\text{resp. }j_{\xi}: \Spec \calO_{X_V, \xi} \longto X_V)
  \end{displaymath}
  be the `inclusion' morphism of the subset $\{\eta\}$ (resp. $\{\xi, \eta\}$) into $X_V$. Let $E \in \calM( \Spec K)$ be an $\calA_K$-module.

  Assume given an extension $\calE \in \calM( \Spec V)$ of $E$. Then $\calE$ embeds canonically into $j_* E$; in particular, the stalk $\calE_{\xi}$ over the discrete valuation ring
  $\calO_{X_V, \xi}$ embeds into the generic fiber $E_{\eta}$. $\calE$ is uniquely determined by $E$ and $\calE_{\xi}$ because
  \begin{equation} \label{intersection}
    \calE = j_* E \cap j_{\xi, *} \calE_{\xi} \subseteq j_{\eta, *} E_{ \eta};
  \end{equation}
  this equation follows easily from the assumption that the special fiber $i^* \calE$ is torsion free, cf. the proof of \cite[Proposition 6]{langton}.

  Moreover, the $\calA$-stable $\calO_{X_V, \xi}$-lattice $\calE_{\xi} \subset E_{\eta}$ is unique up to powers of $\pi$ because its quotient modulo $\pi$ is a simple module under the generic
  fiber of $\calA_l$ by corollary \ref{nonsplit2}. This implies that $\calE$ is determined by $E$ up to isomorphism, thereby proving injectivity.

  To prove surjectivity, we construct an extension $\calE$ of $E$ as follows: The simple $\calA_{\eta}$-module $E_{\eta}$ is Morita equivalent to a one-dimensional vector space over the
  division algebra $D_K = D \otimes_{k( X)} \calO_{X_V, \eta}$. Inside this vector space, we choose a free module of rank one over $D_V := D \otimes_{k( X)} \calO_{ X_V, \xi}$ and denote
  by $\calE_{ \xi}$ the Morita equivalent submodule of $E_{ \eta}$. Then we define $\calE$ by (\ref{intersection}); this clearly defines a sheaf of $\calA_V$-modules over $X_V$ which is flat
  over $V$, whose generic fiber $j^* \calE$ is $E$, and whose special fiber $i^* \calE$ is generically simple. According to the proof of \cite[Proposition 6]{langton} again, $\calE$ is
  coherent over $\calO_{X_V}$, and its special fiber $i^* \calE$ is torsion free. This shows $\calE \in \calM( \Spec V)$.
\end{proof}

\begin{rem} \upshape
  Suppose that $X$ is smooth. In the trivial case $\calA = \calO_X$, generically simple locally projective $\calA$-modules are just line bundles, so they also satisfy the valuative criterion
  for properness; here locally projective and only torsion free modules lie in different connected components of the moduli space. This is no longer true for nontrivial $\calA$, even if
  $\calA$ is a sheaf of Azumaya algebras over an abelian or K3 surface $X$: If $\calA_{\eta}$ is not just a full matrix algebra over $k( X)$, then every generically simple torsion free
  $\calA$-module is a degeneration of locally projective ones by theorem \ref{smooth}.iii below; in particular, the latter do not satisfy the valuative criterion for properness.
\end{rem}

\section{Construction of the coarse moduli schemes} \label{construction}

We choose an ample line bundle $\calO_X(1)$ on $X$ and put $d := \dim( X)$. As usual, the Hilbert polynomial $P( E)$ of a coherent sheaf $E$ on $X$ with respect to this choice is given by
\begin{displaymath}
  P( E; m) := \chi( E( m)) = \sum_{i = 0}^d (-1)^i \dim_k \H^i( X; E( m))
\end{displaymath}
where $E( m) := E \otimes \calO_X(1)^{\otimes m}$. See \cite[Chapter 1.2]{huybrechts-lehn} for details about $P( E)$, in particular for the fact that it is a polynomial
of degree $d = \dim( X)$ if $E$ is torsion free.

Recall that the Hilbert polynomial is locally constant in flat families. Keeping it fixed, we denote by
\begin{displaymath}
  \calM_{\calA/X; P}: \Schemes_k \longto \Sets
\end{displaymath}
the subfunctor of $\calM_{\calA/X}$ that parameterizes families $\calE$ all of whose fibers $\calE_s$ have Hilbert polynomial $P$.

Our first aim is to prove that the class $\calM_{\calA/X; P}( \Spec k)$ of generically simple torsion free $\calA$-modules $E$ with Hilbert polynomial $P$ is bounded
(in the sense of \cite[Definition 1.7.5]{huybrechts-lehn}); this will follow easily from the following apparently weaker statement:
\begin{prop}
  The class of coherent $\calO_X$-modules $E$ with the following two properties is bounded:
  \begin{enumerate}
   \item $E$ is torsion free and has Hilbert polynomial $P$.
   \item $E$ admits at least one $\calA$-module structure for which it is generically simple.
  \end{enumerate}
\end{prop}
\begin{proof}
  Suppose that $E$ is such a sheaf on $X$. We recall a few concepts, which are all in \cite[Chapter 1.2]{huybrechts-lehn}.
  \begin{displaymath}
    \rk( E) := \dim_F( E_{\eta}),\quad \deg( E) := c_1( E) \cdot \calO_X(1)^{d-1}, \quad \quad\mu( E) := \frac{\deg( E)}{\rk( E)}
  \end{displaymath}
  denote the rank, degree and slope of $E$. Writing the Hilbert polynomial as
  \begin{displaymath}
    P( E; m) = \sum_{i = 0}^d \alpha_i( E) \frac{m^i}{i!}
  \end{displaymath}
  with integral coefficients $\alpha_i$ \cite[Chapter 1.2]{huybrechts-lehn}, one has
  \begin{displaymath}
    \deg( E) = \alpha_{d-1}( E) - \rk( E) \alpha_{d-1}( \calO_X),
  \end{displaymath}
  cf. \cite[Definition 1.2.11]{huybrechts-lehn}. As $P( E)$ is fixed, it follows in particular that the slope $\mu( E) = \mu$ is fixed. We denote
  \begin{displaymath}
    \mu_{\max}( E) := \max\{ \mu( E') \big| 0 \neq E' \subseteq E \text{ a coherent $\calO_X$-subsheaf of $E$}\}.
  \end{displaymath}
  $\mu_{\max}( E)$ is in fact the slope of the first subsheaf $E_{\max} \subseteq E$ in the $\mu$-Harder-Narasimhan filtration of $E$ \cite[Section 1.6]{huybrechts-lehn}.

  According to a deep result of A. Langer \cite[Theorem 4.2]{langer}, our class of $\calO_X$-modules $E$ is bounded if the numbers $\mu_{\max}( E)$ are bounded from above. To check the latter,
  we choose an integer $m \in \integers$ such that the coherent sheaf $\calA( m)$ is generated by its global sections. Since $E$ is generically simple, the multiplication map
  \begin{displaymath}
    \calA \otimes_{\calO_X} E_{\max} \longto E
  \end{displaymath}
  is generically surjective. Consequently, the induced map
  \begin{displaymath}
    \H^0( \calA(m)) \otimes_k E_{\max} \longto E( m)
  \end{displaymath}
  is generically surjective, too. But $E$ certainly has a $\mu$-semistable torsion free quotient $E''$ with $\mu(E'') \leq \mu$, e.\,g.\ the last quotient from the $\mu$-Harder-Narasimhan
  filtration of $E$. It is easy to see that there is a nonzero map $E_{\max} \to E''( m)$, obtained by composing
  \begin{displaymath}
    E_{\max} \longto E( m) \longto E''( m).
  \end{displaymath}
  Since $E_{\max}$ and $E''( m)$ are $\mu$-semistable, this implies
  \begin{displaymath}
    \mu_{\max}( E) = \mu( E_{\max}) \leq \mu( E''(m)) = m \mu( \calO_X(1)) + \mu( E'') \leq m \mu( \calO_X(1)) + \mu.
  \end{displaymath}
  This is the required bound for $\mu_{\max}( E)$.
\end{proof}

According to the proposition, there is an integer $m = m_{\calA/X; P}$ with the following property: For every generically simple torsion free $\calA$-module $E$ with Hilbert polynomial $P$,
$E( m)$ is generated by global sections, and $\H^i( E( m)) = 0$ for all $i > 0$. We keep this $m$ fixed in the sequel and denote by $N := P(m)$ the common dimension of all the vector spaces
$\H^0( E(m))$.
\begin{prop} \label{hilb}
 \begin{itemize}
  \item[i)] There is a fine moduli scheme $R$ of finite type over $k$ that parameterizes generically simple torsion free $\calA$-modules $E$ with Hilbert polynomial $P$ together with a basis
   of $\H^0( E(m))$.
  \item[ii)] For $l \gg m$, there is an ample line bundle $L_l$ on $R$ whose fibre at $E$ is canonically isomorphic to $\det \H^0( E( l))$.
  \item[iii)] The algebraic group $\GL( N)$ over $k$ acts on $R$ by changing the chosen bases of the $\H^0( E(m))$.
  \item[iv)] There is a natural action of $\GL( N)$ on $L_l$ that lifts the action in iii.
  \item[v)] The scheme-theoretic stabilizer of every point in $R(k)$ coincides with the scalars $\Gm \subseteq \GL( N)$.
 \end{itemize}
\end{prop}
\begin{proof}
  i) Let $\Quot_P( \calA( -m)^N)$ be Grothendieck's Quot-scheme parameterizing coherent quotients $E$ with Hilbert polynomial $P$ of the $\calO_X$-module sheaf $\calA( -m)^N$. We can take
  for $R$ the locally closed subscheme of $\Quot_P( \calA( -m)^N)$ defined by the following conditions:
  \begin{enumerate}
   \item The quotient sheaf $E$ is torsion free.
   \item $E$ is an $\calA$-module, i.\,e.\ the kernel of $\calA( -m)^N \to E$ is an $\calA$-submodule.
   \item As an $\calA$-module, $E$ is generically simple.
   \item The following composed map is an isomorphism:
    \begin{displaymath}
      k^N \longto \H^0( \calA^N) \longto \H^0( E(m)).
    \end{displaymath}
    Here the left map is given by the unit of the algebra $\calA$.
  \end{enumerate}
  (In particular, this proves that the class of $\calA$-modules $E$ in question is bounded.)

  ii) By Grothendieck's construction of Quot-schemes, the $L_l$ are a fortiori ample line bundles on $\Quot_P( \calA( -m)^N)$.

  iii) and iv) also hold for the whole Quot-scheme, cf. \cite[4.3]{huybrechts-lehn}, and its subscheme $R$ is clearly $\GL( N)$-invariant.

  v) Let $E$ be a generically simple torsion free $\calA$-module with Hilbert polynomial $P$. Choose a basis of $\H^0( E(m))$ and let $G \subseteq \GL( N)$ be the scheme-theoretic
  stabilizer of the corresponding point in $R( k)$. It suffices to show that $G$ and $\Gm$ have the same set of points with values in $k$ and in $k[ \varepsilon]$ where $\varepsilon^2 = 0$.

  Every point in $G(k)$ corresponds to an automorphism of $E$; hence $G( k) = k^*$ because $E$ has only scalar endomorphisms.

  Similarly, every point in $G( k[\varepsilon])$ corresponds to an automorphism of the constant family $E[ \varepsilon]$ over $\Spec k[ \varepsilon]$. Restricting from $\Spec k[ \varepsilon]$
  to $\Spec k$, we get an exact sequence
  \begin{displaymath}
    0 \longto \End_{\calA}( E) \longto[\cdot \varepsilon] \End_{\calA[ \varepsilon]}( E[ \varepsilon] ) \longto \End_{\calA}( E) \longto 0;
  \end{displaymath}
  again because $E$ has only scalar endomorphisms, it implies $\End( E[ \varepsilon] ) = k[ \varepsilon]$ and hence $G( k[ \varepsilon]) = k[ \varepsilon]^*$.

  This proves $G = \Gm$.
\end{proof}
\begin{thm} \label{stable}
  If $l \gg m$, then every point of $R$ is GIT-stable for the action of $\SL( N) \subset \GL( N)$ with respect to the linearization $L_l$.
\end{thm}
\begin{proof}
  We carry the necessary parts of \cite[Chapter 4.4]{huybrechts-lehn} over to our situation.

  Put $V := k^N$; then the points of $R$ correspond to quotients
  \begin{displaymath}
    \rho: V \otimes_k \calA( -m) \to E.
  \end{displaymath}
  We fix such a point.

  Let $V' \subset V$ be a proper vector subspace, and put
  \begin{displaymath}
    E' := \rho( V' \otimes_k \calA( -m)) \subseteq E.
  \end{displaymath}
  Then $E'$ is an $\calA$-submodule of $E$ with nonzero generic fibre; since $E_{\eta}$ contains no proper $\calA_{\eta}$-submodules, this implies $E_{\eta}' = E_{\eta}$.
  So $E$ and $E'$ have the same rank, i.\,e.\ their Hilbert polynomials have the same leading coefficient. Hence
  \begin{equation} \label{estimate}
    \dim( V) \cdot \chi( E'( l)) > \dim( V') \cdot \chi( E( l))
  \end{equation}
  if $l$ is sufficiently large. (We can find one $l$ uniformly for all $V'$ because the family of vector subspaces $V' \subset V$ is bounded.)

  After these preliminaries, we can check that our point $\rho$ in $R$ satisfies the Hilbert-Mumford criterion for GIT-stability, cf.\cite[Theorem 4.2.11]{huybrechts-lehn}. So consider a
  nontrivial one-parameter subgroup $\lambda: \Gm \to \SL( V) = \SL( N)$; we will rather work with the associated eigenspace decomposition $V = \bigoplus_{n \in \integers} V_n$ where $\Gm$
  acts on $V_n$ with weight $n$.

  Let $\bar{R}$ be the closure of $R$ in the projective embedding given by $L_l$; this is also the closure of $R$ in $\Quot_P( V \otimes_k \calA( -m))$ because the $L_l$ are ample on the
  whole Quot-scheme. We have to look at the limit $\lim_{t \to 0} \lambda(t) \cdot \rho$ in $\bar{R}$. This is a fixed point for the $\Gm$-action, so $\Gm$ acts on the fibre of $L_l$ over it,
  necessarily with some weight $- \mu^{L_l}( \rho, \lambda) \in \integers$; what we have to show for stability is $\mu^{L_l}( \rho, \lambda) > 0$.

  First we describe the limit point $\lim_{t \to 0} \lambda(t) \cdot \rho$ as a point in the Quot-scheme. We put
  \begin{displaymath}
    V_{\leq n} := \bigoplus_{\nu \leq n} V_{\nu} \subseteq V \quad\text{and}\quad E_{\leq n} := \rho( V_{\leq n} \otimes_k \calA( -m)) \subseteq E.
  \end{displaymath}
  Then $V_n = V_{\leq n} \big/ V_{\leq n-1}$; we put $E_n := E_{\leq n} \big/ E_{\leq n-1}$, thus obtaining surjections
  \begin{displaymath}
    \rho_n: V_n \otimes_k \calA( -m) \longto E_n.
  \end{displaymath}
  Then
  \begin{displaymath}
    \bar{\rho}:= \bigoplus_{n \in \integers} \rho_n: V \otimes_k \calA( -m) \longto \bigoplus_{n \in \integers} E_n
  \end{displaymath}
  is also a point in $\Quot_P( V \otimes_k \calA( -m))$; it is the limit we are looking for:
  \begin{displaymath}
    \bar{ \rho} = \lim_{t \to 0} \lambda(t) \cdot \rho.
  \end{displaymath}
  To prove this, just copy the proof of \cite[Lemma 4.4.3]{huybrechts-lehn}, replacing $\calO_X( -m)$ by $\calA( -m)$ everywhere.

  The second step is to consider the fibre of $L_l$ over $\bar{\rho}$. It is by definition
  \begin{displaymath}
    \det H^0( \bigoplus_{n \in \integers} E_n( l));
  \end{displaymath}
  this is canonically isomorphic to the tensor product of the determinant of cohomology of the $E_n(l)$. Now $\Gm$ acts on $V_n$ with weight $n$, so it acts on the fiber in question with weight
  \begin{displaymath}
    - \mu^{L_l}( \rho, \lambda) = \sum_{n \in \integers} n \cdot \chi( E_n( l)),
  \end{displaymath}
  cf. \cite[Lemma 4.4.4]{huybrechts-lehn}.

  Finally, we use the preliminaries above to estimate this sum. If we apply (\ref{estimate}) to $V' = V_{\leq n}$, $E' = E_{\leq n}$ and sum up, we get
  \begin{equation} \label{estimate2}
    \sum_{n \in \integers} \big( \dim( V) \cdot \chi( E_{\leq n}( l)) - \dim( V_{\leq n}) \cdot \chi( E( l)) \big) > 0;
  \end{equation}
  note that only finitely many summands are nonzero because $V_{\leq n}$ is zero or $V$ for almost all $n$ since almost all $V_n$ are zero. Put
  \begin{displaymath}
    a_n := \dim( V) \cdot \chi( E_n( l)) - \dim( V_n) \cdot \chi( E( l));
  \end{displaymath}
  again, all but finitely many of these integers are nonzero, and
  \begin{displaymath}
    \sum_{n \in \integers} a_n = \dim( V) \cdot \chi( E( l)) - \dim( V) \cdot \chi( E( l)) = 0.
  \end{displaymath}
  If we write $a_{\leq n} := \sum_{\nu \leq n} a_{\nu}$, then
  \begin{displaymath}
    \sum_n a_{\leq n} + \sum_n n a_n = 0
  \end{displaymath}
  because the sum of the $a_n$ is zero. Hence (\ref{estimate2}) is equivalent to
  \begin{displaymath}
    \sum_{n \in \integers} n \big( \dim( V) \cdot \chi( E_n( l)) - \dim( V_n) \cdot \chi( E( l)) \big) < 0.
  \end{displaymath}
  But $\sum_n n \dim V_n = 0$ because $\Gm$ acts on $V$ with determinant $1$. Thus we obtain
  \begin{displaymath}
    \sum_{n \in \integers} n \big( \dim( V) \cdot \chi( E_n( l)) \big) < 0,
  \end{displaymath}
  i.\,e.\ $- \mu^{L_l}( \rho, \lambda) \cdot \dim( V) < 0$. This proves that the Hilbert-Mumford criterion for GIT-stability is satisfied here.
\end{proof}

Now we can state our main result. See \cite{git} for the concepts `geometric quotient' and `coarse moduli scheme'.
\begin{thm} \begin{itemize}
 \item[i)] The action of $\GL( N)$ on $R$ described in proposition \ref{hilb}.iii above admits a geometric quotient
  \begin{displaymath}
    \coarseM_{\calA/X; P} := \GL( N) \big\backslash R
  \end{displaymath}
  which is a separated scheme of finite type over $k$.
 \item[ii)] The quotient morphism
  \begin{displaymath}
    R \longto \coarseM_{\calA/X; P}
  \end{displaymath}
  is a principal $\PGL( N)$-bundle (locally trivial in the fppf-topology).
 \item[iii)] $\coarseM_{\calA/X; P}$ is a coarse moduli scheme for the moduli functor $\calM_{\calA/X; P}$.
 \item[iv)] $\coarseM_{\calA/X; P}$ is projective over $k$.
\end{itemize} \end{thm}
\begin{proof}
  i) According to geometric invariant theory \cite[Theorem 1.10 and Appendix 1.C]{git} and theorem \ref{stable} above, the action of $\SL( N)$ on $R$ admits a geometric quotient $\coarseM_{
  \calA/X; P}$ which is quasiprojective over $k$. Then $\coarseM_{\calA/X; P}$ is also a geometric quotient for the action of $\GL( N)$ on $R$ because both groups act via $\PGL( N) = \PSL(N)$.

  ii) This morphism is affine by its GIT-construction. According to \cite[Proposition 0.9]{git}, it suffices to show that the action of $\PGL( N)$ on $R$ is free, i.\,e.\ that
  \begin{displaymath}
    \psi: \PGL( N) \times R \longto R \times R, \qquad (g, r) \mapsto (g \cdot r, r)
  \end{displaymath}
  is a closed immersion. The fibers of $\psi$ over $k$-points of $R \times R$ are either empty or isomorphic to $\Spec k$ by proposition \ref{hilb}.v above; furthermore, $\psi$ is proper
  due to \cite[Proposition 0.8]{git}. Hence $\psi$ is indeed a closed immersion. 

  iii) Part i implies that $\coarseM_{\calA/X; P}$ is a coarse moduli scheme for the functor
  \begin{displaymath}
    \underline{\GL( N) \backslash R}: \Schemes_k \longto \Sets
  \end{displaymath}
  that sends a $k$-scheme $S$ to the set of $\GL( N)(S)$-orbits in $R(S)$. However, this functor is very close to the moduli functor $\calM = \calM_{\calA/X; P}$:

  We have a morphism from the functor represented by $R$ to $\calM$; it simply forgets the extra structure. If $S$ is a scheme over $k$, then two $S$-valued points of $R$
  have the same image in $\calM( S)$ if and only if they are in the same $\GL( N)(S)$-orbit. Thus we get a morphism of functors
  \begin{displaymath}
    \phi: \underline{\GL( N) \backslash R} \longto \calM
  \end{displaymath}
  which is injective for every scheme $S$ over $k$. The image of $\phi$ consists of all sheaves $\calE \in \calM( S)$ for which the vector bundle $pr_* \calE( m)$ over $S$ is trivial, where
  $pr: X \times_k S \to S$ is the canonical projection.

  In particular, $\phi$ is bijective whenever $S$ is the spectrum of a field, and it induces an isomorphism between the Zariski sheafifications of both functors. It follows that
  $\coarseM_{\calA/X; P}$ is also a coarse moduli scheme for the functor $\calM$.

  iv) We already know that $\coarseM_{\calA/X; P}$ is quasiprojective. Furthermore, it satisfies the valuative criterion for properness by proposition \ref{proper}.
\end{proof}

In particular, the moduli functor $\calM_{\calA/X}$ of \emph{all} generically simple torsion free $\calA$-modules has a coarse moduli scheme
\begin{displaymath}
  \coarseM_{\calA/X} = \coprod\limits_P \coarseM_{\calA/X; P}
\end{displaymath}
which is a disjoint sum of projective schemes over $k$. If $X$ is smooth of dimension $d$, then we have another such decomposition
\begin{displaymath}
  \coarseM_{\calA/X} = \coprod\limits_{c_1, \ldots, c_d} \coarseM_{\calA/X; c_1, \ldots, c_d}
\end{displaymath}
given by fixing the Chern classes $c_i \in \CH^i( X)$, the Chow group of cycles modulo \emph{algebraic} equivalence. Indeed, each $\coarseM_{\calA/X; c_1, \ldots, c_d}$ is open and closed in
some $\coarseM_{\calA/X; P}$ where $P$ is given by Hirzebruch-Riemann-Roch.

If $X$ is a smooth projective surface, then this decomposition reads
\begin{displaymath}
  \coarseM_{\calA/X} = \coprod\limits_{\genfrac{}{}{0pt}{1}{c_1 \in \NS( X)}{c_2 \in \integers}} \coarseM_{\calA/X; c_1, c_2}.
\end{displaymath}

\section{Deformations and smoothness} \label{deformations}

We introduce the usual cohomology classes that describe deformations of a coherent $\calA$-module $E$, following Artamkin \cite{artamkin1}. By definition, a \emph{deformation} $\calE$ of $E$
over a local artinian $k$-algebra $(A, \m)$ with residue field $k$ is a (flat) family $\calE$ of coherent $\calA$-modules parameterized by $\Spec A$ together with an isomorphism
$k \otimes_A \calE \cong E$.

Consider first the special case $A = k[ \varepsilon]$ with $\varepsilon^2 = 0$. Then we have an exact sequence of $A$-modules
\begin{equation} \label{k_epsilon}
  0 \longto k \longto[\cdot \varepsilon] k[ \varepsilon] \longto k \longto 0.
\end{equation}
By definition, the \emph{Kodaira-Spencer class} of the deformation $\calE$ over $k[ \varepsilon]$ is the Yoneda extension class
\begin{displaymath}
  \ks( \calE) := [0 \longto E \longto[\cdot \varepsilon] \calE \longto E \longto 0] \in \Ext_{\calA}^1( E, E)
\end{displaymath}
obtained by tensoring (\ref{k_epsilon}) over $A$ with $\calE$.
\begin{lemma}
  The Kodaira-Spencer map $\ks$ is a bijection between isomorphism classes of deformations of $E$ over $k[ \varepsilon]$ and elements of $\Ext_{\calA}^1( E, E)$.
\end{lemma}
\begin{proof}
  Let $\calE$ be an $\calA$-module extension of $E$ by $E$. Then $\calE$ becomes an $\calA[ \varepsilon]$-module if we let $\varepsilon$ act via the composition $\calE \twoheadrightarrow E
  \hookrightarrow \calE$. According to the local criterion for flatness \cite[Theorem 6.8]{eisenbud}, $\calE$ is flat over $k[ \varepsilon]$ and hence a deformation of $E$. This defines the
  required inverse map.
\end{proof}

Now let $(A, \m)$ be arbitrary again, and let $\tilde{A}$ be a minimal extension of $A$. In other words, $(\tilde{A}, \tilde{\m})$ is another local artinian $k$-algebra with residue field $k$,
and $A \cong \tilde{A}/(\nu)$ where $\nu \in \tilde{A}$ is annihilated by $\tilde{\m}$. Then we have an exact sequence of $A$-modules
\begin{equation} \label{A_tilde}
  0 \longto k \longto[\cdot \nu] \tilde{\m} \longto A \longto k \longto 0.
\end{equation}
By definition, the \emph{obstruction class} of the deformation $\calE$ over $A$ is the Yoneda extension class
\begin{displaymath}
  \obstr( \calE; k \stackrel{\cdot \nu}{\hookrightarrow} \tilde{A} \twoheadrightarrow A)
    := [0 \longto E \longto[\cdot \nu] \tilde{\m} \otimes_A \calE \longto \calE \longto E \longto 0] \in \Ext_{\calA}^2( E, E)
\end{displaymath}
obtained by tensoring (\ref{A_tilde}) over $A$ with $\calE$. Whenever we want to mention $\calA$, we write $\obstr_{\calA}$ instead of $\obstr$; on the other hand, we may omit
$k \stackrel{\cdot \nu}{\hookrightarrow} \tilde{A} \twoheadrightarrow A$ if they are clear from the context.
\begin{lemma}
  The obstruction class $\obstr( \calE; k \stackrel{\cdot \nu}{\hookrightarrow} \tilde{A} \twoheadrightarrow A)$ vanishes if and only if $\calE$ can be extended to a deformation
  $\tilde{\calE}$ over $\tilde{A}$.
\end{lemma}
\begin{proof}
  If $\tilde{ \calE}$ is a deformation over $\tilde{A}$ extending $\calE$, then we can tensor it over $\tilde{A}$ with the diagram
  \begin{displaymath} \xymatrix{
    0 \ar[r] & k \ar[r]^-{\cdot \nu} \ar@{=}[d] & \tilde{\m} \ar[r] \ar[d] & \m \ar[r] \ar[d] & 0\;\\
    0 \ar[r] & k \ar[r]^-{\cdot \nu}            & \tilde{A}  \ar[r]        & A  \ar[r]        & 0;
  } \end{displaymath}
  this gives us a morphism of short exact sequences of $\calA$-modules
  \begin{displaymath} \xymatrix{
    0 \ar[r] & E \ar[r]^-{\cdot \nu} \ar@{=}[d] & \tilde{\m} \otimes_A \calE \ar[r] \ar[d] & \m \otimes_A \calE \ar[r] \ar[d] & 0\;\\
    0 \ar[r] & E \ar[r]                         & \tilde{ \calE}             \ar[r]        & \calE             \ar[r]         & 0.
  } \end{displaymath}
  The existence of such a morphism implies $\obstr( \calE; k \stackrel{\cdot \nu}{\hookrightarrow} \tilde{A} \twoheadrightarrow A) = 0$ due to the standard exact sequence
  \begin{displaymath}
    \ldots \Ext^1_{\calA}( \calE, E) \longto \Ext^1_{\calA}( \m \otimes_A \calE, E) \longto \Ext^2_{\calA}( E, E) \ldots
  \end{displaymath}

  Conversely, $\obstr( \calE; k \stackrel{\cdot \nu}{\hookrightarrow} \tilde{A} \twoheadrightarrow A) = 0$ implies the existence of such a morphism of short exact sequences of $\calA$-modules.
  Then $\tilde{ \calE}$ becomes an $\tilde{A} \otimes_k \calA$-module if we let any $a \in \tilde{\m} \subset \tilde{A}$ act via the composition $\tilde{\calE} \twoheadrightarrow \calE
  \xrightarrow{a \otimes \_} \tilde{\m} \otimes_A \calE \hookrightarrow \tilde{\calE}$. According to the local criterion for flatness \cite[Theorem 6.8]{eisenbud}, $\tilde{ \calE}$ is flat
  over $\tilde{A}$ and hence a deformation of $E$ extending $\calE$.
\end{proof}

In the special case $A = k[ \varepsilon]$ and $\tilde{A} = k[ \nu]$ with $\nu^3 = 0$, we note that (\ref{A_tilde}) is the Yoneda product of (\ref{k_epsilon}) with itself and hence
\begin{equation} \label{alternating}
  \obstr( \calE; k \stackrel{\cdot \nu}{\hookrightarrow} k[ \nu] \twoheadrightarrow k[ \varepsilon]) = \ks( \calE) \times \ks( \calE).
\end{equation}

From now on, we assume that $X$ is smooth. Then we have a trace map
\begin{displaymath}
  \tr_{\calO_X}: \Ext^i_{\calO_X}( E, E) \longto \H^i( X, \calO_X)
\end{displaymath}
for every coherent $\calO_X$-module $E$: it is defined using a finite locally free resolution of $E$, cf. \cite[10.1.2]{huybrechts-lehn}. If $E$ is a coherent $\calA$-module, then we define
the trace map $\tr = \tr_{\calA/\calO_X}$ as the composition
\begin{displaymath}
  \tr = \tr_{\calA/\calO_X}: \Ext^i_{\calA}( E, E) \longto \Ext^i_{\calO_X}( E, E) \xrightarrow{\tr_{\calO_X}} \H^i( X, \calO_X)
\end{displaymath}
where the first map is induced by the forgetful functor from $\calA$-modules to $\calO_X$-modules. Similarly, one can define a trace map
\begin{displaymath}
  \tr_{\calA/\calO_X}^{\omega_X}: \Ext^i_{\calA}( E, E \otimes_{\calO_X} \omega_X) \longto \H^i( X, \omega_X)
\end{displaymath}
where $\omega_X$ is the canonical line bundle on $X$; cf. \cite[p. 218]{huybrechts-lehn}.

Now suppose that $\calE$ is a deformation of $E$ over $A$. Then we have a line bundle $\det \calE$ over $X_A = X \times_k \Spec A$: it is defined using a finite locally free resolution of
$\calE$ as an $\calO_{X_A}$-module, cf. \cite[1.1.17 and Proposition 2.1.10]{huybrechts-lehn}.
\begin{prop} \label{trace}
  If $\tilde{A}$ is a minimal extension of $A$, then
  \begin{displaymath}
    \tr ( \obstr_{\calA}( \calE; k \hookrightarrow \tilde{A} \twoheadrightarrow A)) = \obstr_{\calO_X}( \det \calE; k \hookrightarrow \tilde{A} \twoheadrightarrow A) \in \H^2( X, \calO_X).
  \end{displaymath}
  In particular, $\tr ( \obstr_{\calA}( \calE)) = 0$ if the Picard variety $\Pic( X)$ is smooth.
\end{prop}
\begin{proof}
  The forgetful map $\Ext_{\calA}^2( E, E) \to \Ext_{\calO_X}^2( E, E)$ maps $\obstr_{\calA}( \calE)$ to $\obstr_{\calO_X}( \calE)$ by definition. It is known that $\tr_{\calO_X}$ maps the
  latter to $\obstr_{\calO_X}( \det \calE)$; cf. Artamkin's paper \cite{artamkin1} for the computation.
\end{proof}

For the rest of this section, we assume that $\calA$ is even a sheaf of Azumaya algebras over the smooth projective variety $X$ of dimension $d$.
\begin{prop} \label{resolution}
  Every coherent sheaf $E$ of $\calA$-modules has a resolution of length $\leq d = \dim( X)$ by locally projective sheaves of $\calA$-modules.
\end{prop}
\begin{proof}
  If $m$ is sufficiently large, then the twist $E( m)$ is generated by its global sections; this gives us a surjection $\partial_0$ of $E_0 := \calA( -m)^N$ onto $E$ for some $N$. Applying the
  same procedure to the kernel of $\partial_0$ and iterating, we obtain an infinite resolution by locally free $\calA$-modules
  \begin{displaymath}
    \ldots E_d \longto[\partial_d] E_{d-1} \ldots E_1 \longto[\partial_1] E_0 \longto[\partial_0] E.
  \end{displaymath}
  We claim that the kernel of $\partial_d$ is locally projective over $\calA$; then we can truncate there, and the proposition follows.

  It suffices to check this claim over the complete local rings $\hat{\calO}_{X, x}$ at the closed points $x$ of $X$; there $\calA$ becomes a matrix algebra $\hat{\calA}_x$, so the resulting
  $\hat{\calA}_x$-modules $\hat{E}_{i, x}$ are Morita equivalent to $\hat{\calO}_{X, x}$-modules. Since $\hat{\calO}_{X, x}$ has homological dimension $d$, the kernel of
  $\partial_{d, x}: \hat{E}_{d, x} \to \hat{E}_{d-1, x}$ is projective over $\hat{\calA}_x$. Hence the kernel of $\partial_d$ is indeed locally projective over $\calA$.
\end{proof}

Our main tool to control the extension classes introduced above will be the following variant of Serre duality. To state it, we fix an isomorphism $\H^d( X, \omega_X) \cong k$.
\begin{prop} \label{serre-duality}
  We still assume that $X$ is smooth of dimension $d$ and that $\calA$ is a sheaf of Azumaya algebras. If $E$ and $E'$ are coherent $\calA$-modules, then the Yoneda product
  \begin{displaymath}
    \Ext_{\calA}^i( E, E') \otimes \Ext_{\calA}^{d-i}( E', E \otimes \omega_X) \longto \Ext^d_{\calA}( E, E \otimes \omega_X)
  \end{displaymath}
  followed by the trace map
  \begin{displaymath}
    \tr_{\calA/\calO}^{\omega_X}: \Ext^d_{\calA}( E, E \otimes \omega_X) \longto \H^d( X, \omega_X) \cong k
  \end{displaymath}
  defines a perfect pairing of finite-dimensional vector spaces over $k$.
\end{prop}
\begin{proof}
  We start with the special case that $E$ and $E'$ are locally projective over $\calA$. Then the $\Ext$-groups in question are Zariski cohomology groups of the locally free $\calO_X$-module
  sheaves $\Homsheaf_{\calA} ( E, E')$ and $\Homsheaf_{\calA}( E', E) \otimes \omega_X$. But $\Homsheaf_{\calA}( E, E')$ and $\Homsheaf_{\calA}( E', E)$ are dual to each other by means of an
  appropriate local trace map, using the fact that the trace map $\calA \otimes_{\calO_X} \calA \to \calO_X$ is nowhere degenerate because $\calA$ is Azumaya. Hence this special case follows
  from the usual Serre duality theorem for locally free $\calO_X$-modules.

  If $E$ and $E'$ are not necessarily locally projective over $\calA$, then we choose finite locally projective resolutions, using proposition \ref{resolution}. Induction on their length
  reduces us to the case where $E$ and $E'$ have resolutions of length one by $\calA$-modules for which the duality in question holds. Now $\Ext^i_{\calA}( E, E')$ and
  $\Ext^{d-i}_{\calA}( E', E \otimes \omega)^{\dual}$ are $\delta$-functors in both variables $E$ and $E'$, and the pairing defines a morphism between them. An application of the five lemma
  to the resulting morphisms of long exact sequences proves the required induction step.
\end{proof}

\begin{thm} \label{smooth}
  Let $X$ be an abelian or K3 surface over $k$, and let $\calA$ be a sheaf of Azumaya algebras over $X$. Suppose $\calA_{\eta} \cong \Mat( n \times n; D)$ for a central division algebra $D$
  of dimension $r^2$ over the function field $k( X)$.
  \begin{enumerate}
   \item[i)] The moduli space $\coarseM_{\calA/X}$ of generically simple torsion free $\calA$-modules $E$ is smooth.
   \item[ii)] There is a nowhere degenerate alternating $2$-form on the tangent bundle of $\coarseM_{\calA/X}$.
   \item[iii)] If $r \geq 2$, then the open locus $\coarseM_{\calA/X}^{\lp}$ of locally projective $\calA$-modules $E$ is dense in $\coarseM_{\calA/X}$.
   \item[iv)] If we fix the Chern classes $c_1 \in \NS( X)$ and $c_2 \in \integers$ of $E$, then
    \begin{displaymath}
      \dim \coarseM_{\calA/X; c_1, c_2} = \Delta/(nr)^2 - c_2( \calA)/n^2 - r^2 \chi( \calO_X) + 2
    \end{displaymath}
    where $\Delta = 2 r^2 n c_2 - (r^2 n- 1) c_1^2$ is the discriminant of $E$.
  \end{enumerate}
\end{thm}
\begin{proof}
  i) We have to check that all obstruction classes
  \begin{displaymath}
    \obstr_{\calA}( \calE; k \hookrightarrow \tilde{A} \twoheadrightarrow A) \in \Ext^2_{\calA}( E, E)
  \end{displaymath}
  vanish. $\Pic( X)$ is known to be smooth; using proposition \ref{trace}, it suffices to show that the trace map
  \begin{displaymath}
    \tr_{\calA/\calO_X}: \Ext^2_{\calA}( E, E) \longto \H^2( X, \calO_X)
  \end{displaymath}
  is injective. But it is straightforward to check that this map is Serre-dual to the natural map
  \begin{displaymath}
    \H^0( X, \omega_X) \longto \Hom_{\calA}( E, E \otimes \omega_X),
  \end{displaymath}
  which is an isomorphism because $\omega_X$ is trivial and $E$ has only scalar endomorphisms.

  ii) Mukai's argument in \cite{mukai} carries over to our situation as follows. We fix an isomorphism $\omega_X \cong \calO_X$. The Kodaira-Spencer map identifies the tangent space
  $T_{[E]} \coarseM_{\calA/X}$ with $\Ext_{\calA}^1( E, E)$. On this vector space, the Serre duality \ref{serre-duality} defines a nondegenerate bilinear form. Indeed, this form is just
  the Yoneda product
  \begin{displaymath}
    \Ext_{\calA}^1( E, E) \otimes \Ext_{\calA}^1( E, E) \longto \Ext_{\calA}^2( E, E);
  \end{displaymath}
  the right hand side is isomorphic to $k$ by Serre duality again. Equation (\ref{alternating}) implies that this bilinear form is alternating because all obstruction classes vanish here.

  iii) Let $E$ be a generically simple torsion free $\calA$-module, and let $\Quot_l( E/\calA)$ be the moduli scheme of quotients $E \twoheadrightarrow T$ where $T$ is a coherent $\calA$-module
  of finite length $l$. This is a closed subscheme of Grothendieck's Quot-scheme $\Quot_{lnr}( E)$ parameterizing those exact sequences of coherent sheaves $0 \to E' \to E \to T \to 0$ for
  which $E'$ is an $\calA$-submodule, i.\,e.\ the composition
  \begin{displaymath}
    \calA \otimes E' \hookrightarrow \calA \otimes E \longto[\cdot] E \twoheadrightarrow T
  \end{displaymath}
  vanishes; here $(nr)^2 = \rk( \calA)$. In particular, $\Quot_l(E/\calA)$ is projective over $k$. We show by induction that $\Quot_l(E/\calA)$ is connected; cf. \cite[6.A.1]{huybrechts-lehn}.

  Let $\Drap_{l_1, l_2}( E/\calA)$ be the moduli scheme of iterated quotients $E \twoheadrightarrow T_1 \twoheadrightarrow T_2$ where $T_i$ is a coherent $\calA$-module of finite length $l_i$
  for $i = 1, 2$; this is again a closed subscheme of some Flag-scheme \cite[2.A.1]{huybrechts-lehn} and hence projective over $k$. Sending such an iterated quotient to $T_1$ and to the pair
  $(T_2, \supp( T_1/T_2))$ defines two morphisms
  \begin{displaymath}
    \Quot_{l+1}( E/\calA) \otgnol[\theta_1] \Drap_{l+1, l}( E/\calA) \longto[\theta_2] \Quot_l( E/\calA) \times X.
  \end{displaymath}
  Using Morita equivalence over the complete local rings at the support of torsion sheaves, it is easy to see that $\theta_1$ and $\theta_2$ are both surjective; moreover, the fibers of
  $\theta_2$ are projective spaces and hence connected. This shows that $\Quot_{l+1}( E/\calA)$ is connected if $\Quot_l( E/\calA)$ is; thus they are all connected.

  Let $E$ still be a generically simple torsion free $\calA$-module; we have to show that its connected component in $\coarseM_{\calA/X}$ contains a locally projective $\calA$-module. Let
  \begin{displaymath}
    E^* := \Homsheaf_{\calO_X}( E, \calO_X)
  \end{displaymath}
  be the dual of $E$; this is a sheaf of right $\calA$-modules. The double dual $E^{**}$ is a sheaf of left $\calA$-modules again; it is locally free over $\calO_X$ and hence locally projective
  over $\calA$. We have an exact sequence
  \begin{displaymath}
    0 \longto E \longto[\iota] E^{**} \longto[\pi] T \longto 0
  \end{displaymath}
  where $T$ is a coherent $\calA$-module of finite length $l$. There is a natural map
  \begin{displaymath}
    \Quot_l( E^{**}/\calA) \longto \coarseM_{\calA/X}
  \end{displaymath}
  that sends a quotient to its kernel; since $\Quot_l( E^{**}/\calA)$ is connected, we may assume that $T$ is as simple as possible, i.\,e.\ that its support consists of $l$ distinct points
  $x_1, \ldots, x_l \in X$ where the stalks $T_{x_i}$ are Morita-equivalent to coherent skyscraper sheaves of length one.

  In this situation, we adapt an argument of Artamkin \cite{artamkin2} to show that $E$ can be deformed to a locally projective $\calA$-module if $r \geq 2$. We consider the diagram
  \begin{displaymath} \xymatrix{
    \Ext^1_{\calA}( E, E) \ar[r]^{\delta} & \Ext^2_{\calA}( T, E) \ar[r]^{\pi^*} \ar[d]^{\iota_*} & \Ext^2_{\calA}( E^{**}, E)\\
    & \Ext^2_{\calA}( T, E^{**}) \ar@{=}[r] & \bigoplus\limits_{i = 1}^l \Ext^2_{\calA}( T_{x_i}, E^{**}).
  } \end{displaymath}
  Here $\pi^*$ is Serre-dual to $\pi^*: \Hom_{\calA}( E, E^{**}) \to \Hom_{\calA}( E, T)$ because $\omega_X \cong \calO_X$. But the only morphisms from $E$ to $E^{**}$ are the multiples of
  $\iota$; hence $\pi^* = 0$, and the connecting homomorphism $\delta$ from the long exact sequence is surjective.

  $\iota_*$ corresponds under Serre duality and Morita equivalence to the direct sum of the restriction maps
  \begin{displaymath}
    \Hom_{\calO_{x_i}}( \calO_{x_i}^r, k_{x_i}) \longto \Hom_{\calO_{x_i}}( \m_{x_i} \oplus \calO_{x_i}^{r-1}, k_{x_i})
  \end{displaymath}
  where $\calO_x = \hat{\calO}_{X, x}$ is the complete local ring of $X$ at $x$, $\m_x \subseteq \calO_x$ is its maximal ideal, and $k_x = \calO_x/\m_x$ is the residue field. Assuming
  $r \geq 2$, these restriction maps are obviously nonzero.

  Hence there is a class $\xi \in \Ext^1_{\calA}( E, E)$ whose image in $\Ext^2_{\calA}( T_{x_i}, E^{**})$ is nonzero for all $i$. Since all obstruction classes vanish, we can find a
  deformation $\calE$ of $E$ over a smooth connected curve whose Kodaira-Spencer class is $\xi$; it remains to show that a general fiber $E'$ of $\calE$ is locally projective over $\calA$.

  Forming the double dual, we get an exact sequence
  \begin{displaymath} 
    0 \longto E' \longto (E')^{**} \longto T' \longto 0.
  \end{displaymath}
  An explicit computation using Morita equivalence shows that the forgetful map
  \begin{displaymath}
    k^r \cong \Ext^2_{\calA}( T_{x_i}, E^{**}) \longto \Ext^2_{\calO_X}( T_{x_i}, E^{**}) \cong k^{n^2 r^3}
  \end{displaymath}
  is injective. Hence the Kodaira-Spencer class of $\calE$ in $\Ext^1_{\calO_X}( E, E)$ also has nonzero image in $\Ext^2_{\calO_X}( T_{x_i}, E^{**})$ for all $i$. According to
  \cite[Corollary 1.3 and the proof of Lemma 6.2]{artamkin2}, this implies that $E'$ is less singular than $E$, i.\,e.\ the length of $T'$ as an $\calO_X$-module is strictly less than $nr$
  at every point of its support. But $T'$ is an $\calA$-module, so these lengths are all divisible by $nr$; hence $T' = 0$, and $E'$ is locally projective over $\calA$.

  iv) Using i and iii, it suffices to compute the dimension of
  \begin{displaymath}
    T_{[E]} M_{\calA/X} \cong \Ext^1_{\calA}( E, E) \cong \H^1( \Endsheaf_{\calA}( E))
  \end{displaymath}
  for a generically simple locally projective $\calA$-module $E$. Note that $\H^0$ and $\H^2$ are Serre-dual to each other and hence both one-dimensional here.

  The endomorphism sheaf $\Endsheaf_{\calA} ( E)$ is an Azumaya algebra of rank $r^2$ over $X$, and the natural map
  \begin{displaymath}
    \calA \otimes_{\calO_X} \Endsheaf_{\calA}( E) \longto \Endsheaf_{\calO_X}( E)
  \end{displaymath}
  is an isomorphism; this is easily checked by reducing to the case that $\calA$ is a matrix algebra. Furthermore, $c_1( \calA)$ is numerically equivalent to zero because
  $\calA \cong \calA^{\dual}$ using the trace over $\calO_X$; similarly, $c_1( \Endsheaf_{\calA}( E))$ and $c_1( \Endsheaf_{\calO_X}( E))$ are also numerically equivalent to zero. Using this,
  the formalism of Chern classes yields
  \begin{equation} \label{c2}
    r^2 c_2( \calA) + (nr)^2 c_2( \Endsheaf_{\calA}( E)) = \Delta
  \end{equation}
  where $\Delta$ is the discriminant of $E$ as above. Hence
  \begin{displaymath}
    \chi( \Endsheaf_{\calA}( E)) = - \Delta/(nr)^2 + c_2( \calA)/n^2 + r^2 \chi( \calO_X)
  \end{displaymath}
  by Hirzebruch-Riemann-Roch.
\end{proof}

\section{Left and right orders} \label{orders}

We still assume that $X$ is smooth projective of dimension $d$ and that $\calA$ is a sheaf of Azumaya algebras over $X$; furthermore, we suppose that the generic fiber $\calA_{\eta}$ is a
division algebra $D$ of dimension $r^2$ over $k( X)$. In this case, generically simple locally projective $\calA$-modules are just locally free $\calA$-modules of rank one; of course every
such $\calA$-module $E$ can be embedded into $D$. The endomorphism sheaf of such a left $\calA$-module $E$ is then just an order in $D$ acting by right multiplication on $E$, exactly as in
the classical picture of the Brandt groupoid \cite[VI, \S2, Satz 14]{deuring}.

The Picard group $\Pic( X)$ acts on the moduli scheme $\coarseM_{\calA/X}$ by tensor product: a line bundle $L \in \Pic( X)$ acts as $E \mapsto E \otimes_{\calO_X} L$. The projective group
scheme $\Pic^0( X)$ of line bundles $L$ algebraically equivalent to zero acts on the individual pieces $\coarseM_{\calA/X; c_1, \ldots, c_d}$ because $c_i( E \otimes L) = c_i( E)$ for all $i$.
The same remarks hold for
\begin{displaymath}
  \coarseM_{\calA/X}^{\lp} = \coprod\limits_{c_1, \ldots, c_d} \coarseM_{\calA/X; c_1, \ldots, c_d}^{\lp}
\end{displaymath}
where the superscript $\lp$ denotes the open locus of locally projective (and hence locally free) $\calA$-modules.
\begin{prop}
  There is a geometric quotient of $\coarseM_{\calA/X}^{\lp}$ by the action of $\Pic( X)$; it is a disjoint sum of separated schemes of finite type over $k$. Its closed points correspond
  bijectively to isomorphism classes of Azumaya algebras over $X$ with generic fiber $D$.
\end{prop}
\begin{proof}
  $\Pic( X)$ acts with finite stabilizers; this follows from
  \begin{equation} \label{dets}
    \det( E \otimes L) \cong \det( E) \otimes L^{\otimes r^2}.
  \end{equation}
  For fixed Chern classes $c_1, \ldots, c_d$, let $G \subseteq \Pic( X)$ be the subgroup of all line bundles $L$ that map $\coarseM_{\calA/X; c_1, \ldots, c_d}^{\lp}$ to
  $\coarseM_{\calA/X; c_1, \ldots, c_d}^{\lp}$. Then $G$ contains $\Pic^0( X)$, and its image in $\NS( X) = \Pic( X)/ \Pic^0( X)$ is contained in the $r^2$-torsion and hence finite,
  so $G$ is a projective group scheme. Therefore a geometric quotient of $\coarseM_{\calA/X; c_1, \ldots, c_d}^{\lp}$ by $G$ exists and is separated and of finite type over $k$,
  according to \cite[Expos{\'e} V, Th{\'e}or{\`e}me 7.1]{sga3}.

  It remains to construct the announced bijection. As
  \begin{displaymath}
    \Endsheaf_{\calA}( E)^{\op} \cong \Endsheaf_{\calA}( E \otimes L)^{\op} =: \calA',
  \end{displaymath}
  which is again an Azumaya algebra with generic fiber $D$, we obtain a well defined map from closed points of the quotient to isomorphism classes of such $\calA'$. Conversely, given
  $\calA, \calA'$, the possible locally free $\calA$-modules $E$ of rank one with $\Endsheaf_{\calA}( E)^{\op} \cong \calA'$ all differ only by tensoring with line bundles. This can be seen
  as follows:

  Suppose that $E$ and $E'$ are locally free $\calA$-modules of rank one with
  \begin{displaymath}
    \Endsheaf_{\calA}( E)^{\op} \cong \Endsheaf_{\calA}( E')^{\op} \cong \calA'.
  \end{displaymath}
  We choose embeddings of $E$ and $E'$ into $\calA_{\eta} = D$; this also embeds $\Endsheaf_{\calA}( E)^{\op}$ and $\Endsheaf_{\calA}( E')^{\op}$ into $D$. The given isomorphism between them
  induces an automorphism of $D$, i.\,e.\ conjugation with an element of $D$; altering the embedding $E' \hookrightarrow D$ by a right multiplication with this element, we may assume that
  $\Endsheaf_{\calA}( E)^{\op} = \Endsheaf_{\calA}( E')^{\op} =: \calA'$ as subalgebras of $D$.

  There is an open subscheme $U \subseteq X$ such that $\calA|_U = E|_U = E'|_U = \calA'|_U$. Furthermore, $X \setminus U$ is a finite union of divisors $D_1, \ldots D_l$; it is enough to
  study the question at the generic points $x_i$ of the $D_i$. There the local ring $\calO_{X, x_i}$ is a discrete valuation ring; over its completion $\hat{\calO}_{X, x_i}$, $\calA$ becomes
  isomorphic to a matrix algebra, so we can describe the situation using Morita equivalence as follows:

  $\hat{E}_{x_i}$, $\hat{E}_{x_i}'$ correspond to lattices over $\hat{\calO}_{X, x_i}$ in $F_{x_i}^r$ ($F_{x_i}$ the completion of $F$ at $x_i$) such that
  $\End_{\hat{\calO}_{x_i}}( \hat{E}_{x_i})^{\op} = \End_{\hat{\calO}_{x_i}}( \hat{E}_{x_i}')^{\op} = \hat{ \calA}_{x_i}'$. But then it is an easy exercise to see that
  $\hat{E}_{x_i}' = \pi_i^{N_i} \hat{E}_{x_i}$ for some $N_i \in \integers$ where $\pi_i \in \hat{\calO}_{x_i}$ is a uniformising element. From this our claim follows.

  This shows that the map above is injective. For the surjectivity, we consider two Azumaya algebras $\calA, \calA' \subseteq D$ and define
  $E(U):= \{f \in D: \calA|_U \cdot f \subseteq \calA'|_U\}$. Using Morita equivalence as above, it is easy to check that $E$ is a locally free $\calA$-module of rank one
  with $\Endsheaf_{\calA}( E)^{\op} = \calA'$; this proves the surjectivity.
\end{proof}

\begin{rem} \upshape
  If $\calA'$ is another sheaf of Azumaya algebras over $X$ with generic fiber $\calA_{\eta}' \cong D = \calA_{\eta}$, then the moduli spaces $\coarseM_{\calA/X}$ and $\coarseM_{\calA'/X}$
  are isomorphic. Indeed, the preceeding proof shows that there is a locally free $\calA$-module $E$ of rank one with $\Endsheaf_{\calA}( E)^{\op} \cong \calA'$; then $E$ is a right
  $\calA'$-module, and one checks easily that the functor $E \otimes_{\calA'} \_$ defines an equivalence from left $\calA'$-modules to left $\calA$-modules.
\end{rem}

\begin{rem} \upshape
  If $X$ is a surface, then this quotient can be decomposed explicitly into pieces of finite type as follows:

  The action of $\Pic( X)$ preserves the discriminant $\Delta( E) \in \integers$, so we get a decomposition
  \begin{displaymath}
    \coarseM_{\calA/X}^{\lp} \big/ \Pic( X) = \coprod_{\Delta \in \integers} \coarseM_{\calA/X; \Delta}^{\lp} \big/ \Pic( X).
  \end{displaymath}
  Now the first Chern class $c_1( E) \in NS( X)$ decomposes $\coarseM_{\calA/X; \Delta}^{\lp}$ into pieces of finite type. But $c_1( E \otimes L) = c_1( E) + r^2 c_1(L)$, and $r^2 \NS( X)$
  has finite index in $\NS( X)$, so $\coarseM_{\calA/X; \Delta}^{\lp} \big/ \Pic( X)$ is indeed of finite type over $k$.

  According to equation (\ref{c2}), fixing $\Delta( E)$ corresponds to fixing $c_2( \calA') \in \integers$ where $\calA' = \Endsheaf_{\calA}( E)^{\op}$. If $X$ is an abelian or K3 surface,
  then theorem \ref{smooth}.iv yields
  \begin{displaymath}
    \dim \coarseM_{\calA/X; \Delta}^{\lp} \big/ \Pic( X) = c_2 - (r^2 -1) \chi( \calO_X)
  \end{displaymath}
  where $c_2 = c_2( \calA') = \Delta/r^2 - c_2( \calA) \in \integers$ is the second Chern class of the Azumaya algebras $\calA'$ that this quotient parameterizes.
\end{rem}

\begin{rem} \upshape
  M. Lieblich \cite{lieblich1, lieblich2} has compactified such moduli spaces of Azumaya algebras using his generalized Azumaya algebras, i.\,e.\ algebra objects in a derived category
  corresponding to endomorphism algebras of torsion free rank one $\calA$-modules.
\end{rem}

\end{document}